\documentstyle[11pt,amsfonts,amssymb,twoside]{article}

\newtheorem{th}{Theorem}
\newtheorem{Prop}{Proposition}[section]
\newtheorem{Prop(nosectionno)}{Proposition}
\newtheorem{lemma}{Lemma}
\newenvironment{co}{\begin{trivlist}\item[]{\bf Corollary\ }\begin{em}}
{\end{em}\end{trivlist}}

\newenvironment{pf}{\begin{trivlist}\item[]{\bf Proof:}}
{\rule{.08in}{.08in}\end{trivlist}}

\begin{document}
\title{Rank One symmetric Spaces and Rigidity}
\author{Inkang Kim\\University of California at Berkeley}
\date{}
\maketitle
\begin{abstract}
In this paper we show that if we have two nonelementary,irreducible
representations from some group G(not necessarily finitely generated) into
$Iso^{\pm}(H^n_F)$ and $Iso^{\pm}(H^m_K)$ such that $dim_R(F) \times n=dim_R(K)
\times m $ where $F,K=R,C,H,O$ and have the same marked length spectrum
,then $F=k,n=m$ and two representations are conjugate.In the last section we
show that if we have  a representation from a finitely generated group into
$SL_2(C)$,sitting smoothly in the space of representations then we can find
smooth coordinate chart by translation lengthes of finite number of hyperbolic
elements around this representation.

\end{abstract}
\openup6pt
\section{Introduction}
Let F to be $R,C,H$.The inner product on $F^{n+1}$ is given as;
$$<(z_1,\cdots,z_{n+1}),(w_1,\cdots,w_{n+1})>=
\sum_1^n{z_i\bar{w_i}}-z_{n+1}\overline{w_{n+1}}$$
Let $GL(n+1,F)$ act on $F^{n+1}$ from the right and let $P^n_F$ be the space of
left $F$-lines.
Then $G=O_F(n,1)$ is a subgroup of $GL(n+1,F)$ preserving this inner product.
The induced action of $O_F(n,1)$ on $P^n_F$ gives three invariant subsets
$D_+,D
_0,D_-$ and the action is positive definite on $D_+$.A component of $D_+$ is
called hyperbolic space.An isotropy group of $O_F(n,1)$ is a maximal compact
 subgroup
 $K=O_F(1) \times O_F(n)$.If we numerate them;
    $$ O_R(n,1)=O(n,1),K=O(1)\times O(n)$$
    $$     O_C(n,1)=U(n,1),K=U(1)\times U(n)$$
   $$  O_H(n,1)=Sp(n,1),K=Sp(1) \times Sp(n)$$

The hyperbolic Cayley plane is similarly realized and the isometry group is the
real form of $F_4$ of real rank one.The maximal compact subgroup is $Spin(9)$.

We can compactify the hyperbolic space by adding $D_0$.If $dim_RF=k$ then $D_0$
is a $(kn-1)$ dimensional sphere denoted by $\partial{H^n_F}$.

Take the Iwasawa decomposition $G=KAN$.Then $A$ and $N$ fix $\xi$ in the sphere
at infinity and $N$ is a nilpotent goup obtained as an extion
 $$ 0 \rightarrow  ImF \rightarrow N \rightarrow F^{n-1} \rightarrow 0.$$
Horospheres based at $\xi$ are the orbits of $N$.
The stabilizer in $O_F(n,1)$ of $\xi$ is a parabolic subgroup with Langlands
decomposition $MAN$.For $F=R,C,H$ $M$ is a subgroup of $K$ of the form
$O_F(1)\times O_F(n-1)$.For Cayley plane $M$ is $Spin(7)$.
$K$ actually acts transitively on $\partial{H^n_F}$ giving
$K/M=\partial{H^n_F}$
.The centralizer of $M$ in $K$ acts freely on $K/M$ with quotient $P^{n-1}_F$.
This gives Hopf fibration $\partial{H^n_F}\rightarrow P^{n-1}_F$.

Since $K$ acts transitively on $\partial{H^n_F}$ ,killing form on its Lie
algebra determines a Riemannian metric on $\partial{H^n_F}$.Define the
orthogonal complements
 to the vertical tangent spaces of the Hopf fibration $\partial{H^n_F}
\rightarrow P^{n-1}_F$.Then there is $K$ invariant Riemannian metric on the
horizontal bundle  and if we define the distance of two points
 as the minimum of the lengthes of the pathes joining these two points and
staying horizontal bundle it is Carnot-Caratheodary metric.Mitchell[Mi]
 calculated
that the housdorff dimension of $\partial{H^n_F}$ in this metric is $k(n+1)-2$.

\section{Boundary of Rank one symmetric Spaces}
The boundary of rank one symmetric space  $H^m_F$ is a one point
compactification of nilpotent group $N$ in Iwasawa decomposition,denoted by $N
\cup \infty$.We will introduce nice coordinates and left invariant distance to
introduce cross ratio of four points.We will define everything in Cayley
numbers because that is the most
general case among four.

Cayley number is a pair of quaternions $(q_1,q_2)$.The multiplication is
defined
 as:
$$(q_1,q_2)(p_1,p_2)=(q_1p_1-\bar{p_2}q_2,q_2p_1+p_2\bar{q_1})$$

Also we define $\overline{(q_1,q_2)}=(\bar{q_1},-q_2)$.

Then it has the following properties.

$1)q\bar{q}=|q|^2$

$2)|qp|=|q||p|$

$3)q^{-1}=\bar{q}/|q|^2$

$4)\overline{pq}=\bar{q}\bar{p}$

Even though Cayley numbers are not commutative nor associative,using Artin's
theorem saying that the subalgebra generated by two elements
is associative,we have following extra properties.

1)$(xy)y^{-1}=x$

2)$(xy)^{-1}=y^{-1}x^{-1}$

Here is some usuful fact. $Aut(O)=G_2$ and it is compact group.It fixes real
number and acts transitively on unit pure imaginaries.The stabilizer of $i$ is
a
 copy of $SU(3)$ and it acts transitively on the unit pure imaginaries
orthogona
l to $i$.The stabilizer of $i,j$ which fixes $k$ since $k=ij$,acts transitively
on the unit pure imaginaries orthogonal to $i,j,k$ and it is a copy of $SU(2)$.
For more extensive informations,see [Fr].

Now an element of $N$ is denoted as $[(t,q),k]$ where

 1)$F=R(Real),(t,q)=0,k \in R^{m-1}$.

2)$F=C(Complex),t$ is real,$q=0,k \in C^{m-1}$.

3)$F=H(Quaternion),t \in R^3,q=0,k \in H^{m-1}$.

4)$F=O(Cayley),t \in R^3,q$ is a quaternion,k is a Cayley number.

The multiplication is defined as;
$$[(t,q),k][(t',q'),k']=[(t+t',q+q')+2Im<k,k'>,k+k']$$

We define the gauge of $[(t,q),k]$ as;
 $$A([(t,q),k])=(|k|^2+t,q)$$

Then define $|[(t,q),k]|=(|k|^4+|t|^2+|q|^2)^{1/4}$ and $d(g,g')=|g'^{-1}g|$.
This is a left invariant distance and the housdorff dimesion of the boundary
with respect to this metric is $dim_R(F) \times (m+1)-2$ wich agrees with
Mitchell's calculation.Even though this metric is not riemannian if we take the
inner
metric of this it becomes Carnot-Caratheodory metric.For references
see[Gr][Pa][Mi].

The action of $Iso(H^m_F)$ extends continously to $\partial(H^m_F)$.Let
$Sim(N)$ denote the subgroup of $Iso(H^m_F)$ wich fixes $\infty$.Then it is
isomorphic to $ N \rtimes (O_F(m-1)\cdot O_F(1) \times R)$ and $N \rtimes
O_F(m-1)\cdot O_F(1)$ acts as  isometries with respect to the metric given
above.Note that $N$ part comes from left action of the group itself and
$O_F(m-1)\cdot O_F(1) \times R$ part comes from
hyperbolic isometries of the form in $O_F(m,1)$
\[
\left[
\begin{array}{ccc}
   M   &   0   &  0  \\
   0   &   \nu\cosh(s)  & \nu\sinh(s)  \\
   0  &    \nu\sinh(s)  & \nu\cosh(s)
\end{array}
\right]
\]
where $O_F(m-1)\cdot O_F(1)=O(m-1),U(m-1),Sp(m-1)\cdot Sp(1),O_O(1)\cdot
O_O(1)$ depending on $F=R,C,H,O$ respectively.These actions will be described
in next
section.

We define the cross ratio of four points as;
$$[g_1,g_2,g_3,g_4]=\frac{|A(g_3^{-1}g_1)||A(g_4^{-1}g_2)|}
{|A(g_4^{-1}g_1)||A(g_3^{-1}g_2)|}$$

Let's diverge to the unit ball model for a while so that we can get some
connetion between the metric of hyperbolic space and the metric defined above.
In Real,complex,Quaternion case for $x,y \in H^m_F$
$$cosh(d(x,y))=\frac{|1-<x,y>|}{(1-<x,x>)^{1/2}(1-<y,y>)^{1/2}}$$

For Cayley hyperbolic case if we define
$$R<v,w>=Re(v_1\bar{v_2})(w_2\bar{w_1})-Re(\bar{v_2}w_2)(\bar{w_1}v_1)$$

Then for $x,y \in H^2_O$ the distance is
$$cosh(d(x,y))=\frac{(|1-<x,y>|^2+2R<x,y>)^{1/2}}
{(1-<x,x>)^{1/2}(1-<y,y>)^{1/2}}$$
For details see[Mo].
Let $<<x,y>>=1-<x,y>$ for $F=R,C,H$ and $<<x,y>>={(|1-<x,y>|^2+2R<x,y>)^{1/2}}$
for Cayley case.Then the cross ratio of the four points in the boudary of
hyperbolic space is defined as
 $$[x,y,z,w]=\frac{<<z,x>><<w,y>>}{<<w,x>><<z,y>>}$$
Note that this is a limit of
$$\frac{\cosh(d(z_i,x_i))\cosh(d(w_i,y_i))}
{\cosh(d(w_i,x_i))\cosh(d(z_i,y_i))}$$
where $x_i,y_i,z_i,w_i$ tend to $x,y,z,w$.So it is invariant under isometries
of hyperbolic space.

To make connection between these two definitions of cross ratio we introduce
generalize projection from $N \cup \infty$ to the boundary of unit ball.
$$w_1=2(1+|k|^2-t,-q)^{-1}k$$
$$w_2=(1+|k|^2-t,-q)^{-1}(1-|k|^2+t,q)$$
Note that $(0,0)$ corresponds to $(0,1)$,$\infty$ to $(0,-1)$.

Let $g_i=[(t_i,q_i),(c_i,d_i)]$ for $i=1,2$.
By the general projection these two points correspond to
$$x_i=1/((1+|k_i|^2)^2+|t_i|^2+|q_i|^2)
[2(c_i+|k_i|^2c_i+t_ic_i-\bar{d_i}q_i,q_ic_i+d_i+$$
  $$d_i|k_i|^2-d_it_i),(1-|k_i|^4+2t_i-|t_i|^2-|q_i|^2,2q_i)]$$

Since (0,0) and $\infty$ correspond to $(0,1)$ and $(0,-1)$
$$<<x_i,(0,-1)>>=2/((1+|k_i|^2)^2+|t_i|^2+|q_i|^2)$$
$$<<x_i,(0,1)>>=2\frac{\sqrt{(|k_i|^4+|k_i|^2+|t_i|^2+|q_i|^2)^2
+|t_i|^2+|q_i|^2}}{(1+|k_i|^2)^2+|t_i|^2+|q_i|^2}$$

Since
$((|k_i|^4+|k_i|^2+|t_i|^2+|q_i|^2)^2+|t_i|^2+|q_i|^2=(|k_i|^4+|t_i|^2+|q_
i|^2)((1+|k_i|^2)^2+|t_i|^2+|q_i|^2)$ we get

$$[\infty,0,x_1,x_2]=\frac{(k_2|^4+|t_2|^2+|q_2|^2)^{1/2}}
{(k_1|^4+|t_1|^2+|q_1|^2)^{1/2}}$$
$$=|A(g_2^{-1})|/|A(g_1^{-1})|$$
so we are done.

For the case $[0,g_1,\infty,g_2]$,just note that
$[0,g_1,\infty,g_2]=[g_2^{-1},g_2^{-1}g_1,\infty,0]$.By using above fact we are
done again.

In fact for $F=R,C,H$ it is easy to show that two definitions of cross ratio
agree by brutal calculation but in Cayley Hyperbolic case calculation is quite
hard.But to prove our main theorem we will just need what we proved.
\section{Action of Isometries on $\partial{H^m_F}$}

In this section we want to incode the action of hyperbolic isometry on the
boundary.A hyperbolic isometry in $O_F(m,1)$ has the form of ;
\[
\left[
\begin{array}{ccc}
   M   &   0   &  0  \\
   0   &   \nu\cosh(s)  & \nu\sinh(s)  \\
   0  &    \nu\sinh(s)  & \nu\cosh(s)
\end{array}
\right]
\]
where $\nu=1,M \in O_F(m-1)$ for $F=R,C$ and $|\nu|=1,M \in Sp(m-1)$ for $F=H$
and $ |M|=|\nu|=1$ for Cayley hyperbolic case.

Then this element send $[(t,q),k]$ to;
$$(t',q')=e^{-2s}{\nu}^{-1}(t,q)\nu$$
\begin{eqnarray*}
k'=e^{-s}\{[{\nu}^{-1}(e^{2s}+|k|^2-t,-q)^{-1}{\nu}]^{-1}\}& &  \\
\{[{\nu}^{-1}((e^{2s}+|k|^2-t,-q)^{-1}(1+|k|^2-t,-q))][((1+|k|^2-t,-q)^{-1}k)A]
\}& &
\end{eqnarray*}

First note that the matrix acts on the unit ball as;
$$w_1'=(w_2{\nu}\sinh(s)+{\nu}\cosh(s))^{-1}(w_1A)$$
$$w_2'=(w_2{\nu}\sinh(s)+{\nu}\cosh(s))^{-1}(w_2\nu\cosh(s)+\nu\sinh(s))$$
using coordinate chage between  unit ball model and $F^{m+1}$ sending $w$ to
$(w,1)$.
To check above claim we use Artin's theorem crucially.
By using generalized projection and above equations we can show
\begin{eqnarray*}
w_2'=\{((1+|k|^2-t,-q)^{-1}(1-|k|^2+t,q))\nu\sinh(s)+\nu\cosh(s)\}^{-1}& & \\
 \{[(1+|k|^2-t,-q)^{-1}(1-|k|^2+t,q)]\nu\cosh(s)+\nu\sinh(s)\}& &
\end{eqnarray*}
$$=[{\nu}^{-1}((e^{2s}+|k|^2-t,-q)^{-1}(1+|k|^2-t,-q))]
[((1+|k|^2-t,-q)^{-1}(e^{2s}-|k|^2+t,q))\nu]$$  equal to
$$[({\nu}^{-1}(e^{2s}+|k|^2-(t,q))^{-1})\nu]
[{\nu}^{-1}(e^{2s}-|k|^2+(t,q)\nu]$$

To prove this ,set $r_1=e^{2s}+|k|^2,r_2=e^{2s}-|k|^2,r_3=1+|k|^2,(t,q)=Q$ and
using identity $-QQ=|Q|^2$
we can show
$$[{\nu}^{-1}((r_1+Q)^{-1}(r_3-Q))][((r_3+Q)(r_2+Q))\nu]=
[({\nu}^{-1}(r_1-Q)^{-1})\nu][{\nu}^{-1}(r_2+Q)\nu].$$
For $w'_1$,it is easy to check.
\section{Marked length spectrum determines representation}
First we will make very simple but important observation to prove the theorem.
\begin{lemma}Let $a,b$ be two hyperbolic isometries such that $x_1$ is a
repelling
fixed point of $a$ ,$x_3$ is a attracting fixed point of $a$ and $x_4$ is a
attracting fixed point of $b$,$x_2$ is a repelling fixed point of $b$.Then
 $$\lim_{n\rightarrow\infty}e^{l(a^n)+l(b^n)-l(a^nb^n)}=|[x_1,x_2,x_3,x_4]|.$$
\end{lemma}
\begin{pf}Choose $x_1^n$ on the axis of $a$,$z_1^n$ on the axis of
$b^na^n$,$x_2
^n$ on the axis of $b$,$w_2^n$ on the axis of $a^nb^n$ so that $d(x_1^n,z_1^n)$
and $d(x_2^n,w_2^n)$ go to zero.If we put
$x_3^n=a^n(x_1^n),z_3^n=a^n(z_1^n),x_4
^n=b^n(x_2^n),w_4^n=b^n(w_2^n)$ and
$d_{13}^n=d(x_1^n,x_3^n),d_{24}^n=d(x_2^n,x_
4^n),d_{23}^n=d(w_2^n,z_3^n),
d_{14}^n=d(z_1^n,w_4^n)$ then
$$\lim \sqrt{\frac{<<x_1^n,x_3^n>><<x_3^n,x_1^n>>
<<x_2^n,x_4^n>><<x_4^n,x_2^n>>}
{<<x_1^n,x_4^n>><<x_4^n,x_1^n>><<x_2^n,x_3^n>><<x_3^n,x_2^n>>}}$$
$$=\lim
\frac{(e^{d^n_{13}}+e^{-d^n_{13}})(e^{d^n_{24}}+e^{-d^n_{24}})}{(e^{d^n_
{14}}+e^{-d^n_{14}})(e^{d^n_{23}}+e^{-d^n_{23}})}$$
$$=\lim \frac{e^{d^n_{13}}e^{d^n_{24}}}{e^{d^n_{14}}e^{d^n_{23}}}+\lim
\frac{e^{
d^n_{13}}}{e^{d^n_{14}}e^{d^n_{23}}}+\lim
\frac{e^{d^n_{24}}}{e^{d^n_{14}}e^{d^n
_{23}}}$$
$$=\lim e^{d^n_{13}+d^n_{24}-d^n_{14}-d^n_{23}}+\lim
e^{d^n_{13}-d^n_{14}-d^n_{2
3}}+\lim e^{d^n_{24}-d^n_{14}-d^n_{23}}$$
$$=\lim e^{l(a^n)+l(b^n)-l(a^nb^n)}+\lim e^{l(a^n)-l(a^nb^n)}+\lim
e^{l(b^n)-l(a
^nb^n)}$$
Since $l(a^n)+l(b^n)-l(a^nb^n)$ always
exists,$l(a^n)-l(a^nb^n),l(b^n)-l(a^nb^n)
$ both go to $-\infty.$ For more general argument see [Kim].

So we finally get
$$\lim e^{l(a^n)+l(b^n)-l(a^nb^n)}=[x_1,x_2,x_3,x_4]$$

\end{pf}

This lemma shows that marked length spectrum determines cross ratio on the
limit set because every two points in the limit set can be approximated by end
points of some hyperbolic isometry.By using limit argument it is clear that
cross ratio of every four points in the limit set is determined by its marked
length spectrum.
Now we are ready to prove the theorem.
\begin{th}Let $\rho:G\rightarrow Iso^{\pm}(H^n_K),\phi:G\rightarrow
Iso^{\pm}(H^
m_F)$ be two nonelementary,irreducible representations having the same marked
le
ngth spectrum such that $n\times dim_R(K)=m \times dim_R(F)$ where $K,F$ is
real
,complex,quaternion or cayley fields.Then $n=m,K=F$ and they are conjugate.
\end{th}
\begin{pf}
After we conjugate representations we may assume that 0,and $\infty$ are in the
limit sets and they are the two fixed points  of hyperbolic
 isometries $a$ and $a'$ in $\rho$ and $\phi$ respectively.Furthermore we have
$$\frac{|x_i|}{|x_j|}=|[x_i,x_j,[0,0],\infty]|^{1/2}=|[y_i,y_j,[0,0],\infty]|^{
1/2}=\frac{|y_i|}{|y_j|}$$
So after scaling we can assume that $|x_i|=|y_i|$.
Similarly
$$\frac{d(x_i,x_j)}{|x_j|}=|[[0,0],x_i,\infty,x_j]|^{1/2}=|[[0,0],y_i,
\infty,y_j
]|^{1/2}=\frac{d(y_i,y_j)}{|y_j|}$$
and hence $d(x_i,x_j)=d(y_i,y_j).$

Let [0,1],$[0,z]$ are the two points in the limit set of $\rho$ and
$[c,d],[(t,q
),k]$ be the corresponding points in the limit set of $\phi$.

Then $a$ sends $[(t,q),k]$ to
$$(t',q')=e^{-2s}{\nu}^{-1}(t,q)\nu$$
\begin{eqnarray*}
k'=e^{-s}\{[{\nu}^{-1}(e^{2s}+|k|^2-t,-q)^{-1}{\nu}]^{-1}\}& &  \\
\{[{\nu}^{-1}((e^{2s}+|k|^2-t,-q)^{-1}(1+|k|^2-t,-q))][((1+|k|^2-t,-q)^{-1}k)M]
\}& &
\end{eqnarray*}
where $\nu=1,M \in O(m-1),U(m-1)$ in real and complex case,$|\nu|=1,M \in
Sp(m-1
),|M|=1$ in quaternion and cayley cases.

Similarly $a'$ sends $[(t,q),k]$ to
$$(t',q')=e^{-2s}{\nu}^{-1}(t,q)\nu$$
\begin{eqnarray*}
k'=e^{-s}\{[{\mu}^{-1}(e^{2s}+|k|^2-t,-q)^{-1}{\mu}]^{-1}\}& &  \\
\{[{\mu}^{-1}((e^{2s}+|k|^2-t,-q)^{-1}(1+|k|^2-t,-q))][((1+|k|^2-t,-q)^{-1}k)N]
\}& &
\end{eqnarray*}

Now
 $$d([0,z],[0,1])^4=d([(t,q),k],[c,d])^4$$ gives us that $|(t,q)|$ is a
function
 of $z,k,c,d$.But
$$d([0,z],a^n([0,1]))^4=d([(t,q),k],a'^n([c,d]))^4$$ gives that $|(t,q)|$ is a
function of $e^{-ns},z,k,c,d$ for all n.This is because right side of equation
has the higher power of $e^{-s}$.This is not possible unless $|(t,q)|=0$.

By this way we can conclude that all the limit sets of $\phi$ corresponding to
$[0,z]$ of $\rho$ are the form of $[0,k]$.Then by further conjugation we can
assume that [0,1] is also the corresponding point of $\phi$ to [0,1] of $\rho$.

Next we want to show that either $M=N,\nu=\mu$ or $M=\bar{N},\nu=\bar{\mu}$.
This can be seen from the equation
$$d(a^n([0,1]),[0,1])^4=d(a'^n([0,1]),[0,1])^4$$ for all n.
Now by conjugating the representation by the map $(k \rightarrow \bar{k})$ we
can assume that $M=N,\nu=\mu$.

If $[0,w]$ is a point in the limit set of $\rho$ and $[0,w']$ is the
corresponding point of $\phi$ then
$$d(a^n([0,w]),[0,1])^4=d(a'^n([0,w']),[0,1])^4$$
gives that $w=w'$.

These all show that $K=F,m=n$.So far we conjugated two representations to get
that two limit sets are equal on the set $K^{m-1}=\{[0,k]:k \in
R^{m-1},C^{m-1},Q^{m-1},
O\}$ and two spaces are actually the same one of the four hyperbolic spaces.

For any point $x$ of the limit set of $\rho$      not in the set $K^{m-1}$
 there are only two possible corresponding points $y$(either $x$ itself or
reflection of $x$ along $K^{m-1}$)
 such that the distance between $x$ and every point in the limit set of $\rho$
  lying on the set $K^{m-1}$ is equal to the
distance between $y$ and
corresponding point of $\phi$.
By taking reflection along this set if necessary,we can assume that $x=y$.
Now for any point in the limit set of $\rho$ off the set $K^{m-1} \cup x$ there
is a unique corresponding point as above.So two limit sets are actually equal.
This shows that we found an isometry such that after we conjugate one
representation by this isometry the limit sets are equal.But this actually
shows that two representations are equal.The reason is the following.
Note that the end points of hyperbolic isometry $aba^{-1}$ are the images of
end
 points of $b$ under $a$.But since every isometry is determined by  images of
finitely many points on $\partial{H^m_F}$ and since two representations have
the same limit sets after conjugation each element should be the same.
\end{pf}
\begin{co}Let M,N be rank one locally symmetric manifolds of the same dimension
and homotopically equivalent.If none of them have totally geodesically embedded
submanifold of dimension greater than one and if they have the same marked
length spectrum then they are isometric.
\end{co}

Since three manifolds get special attention we will write down the corollary
for hyperbolic three manifolds.Here we do not need the assumption about
irreducibility.
\begin{co}If we have two nonelementary orientable hyperbolic three manifolds
having the same marked length spectrum then they are isometric.
\end{co}
\begin{pf}Since every orientation preserving Mobius map is determined by the
images of three distinct points whether they are on the line or not we don't
have to worry about irreducibility.
\end{pf}
\section{Space of Representations into $SL_2(C)$}
\subsection{Backgroud}
Let G be a finitely generated group.Space of representations from G to
$SL_2(C)$
 is a set of homomorphisms from G into $SL_2(C)$.We say that two
representations
 are equivalent if they are conjugated by some element of $SL_2(C)$.The
character of  a representation $\rho$ is the function $\chi_{\rho}:G
\rightarrow C$ such
 that $\chi_{\rho}(g)=tr(\rho(g))$.
If $G=\{g_i,i=1,\cdots,n;r_1=\cdots=r_k=1\}$  then space  of representations
$R(G)$ is
a subset of $SL_2(C) \subset C^{4n}$,wich are set of all points $(\rho(g_1),
\cdots,\rho(g_n))$ satisfying the relator relations.It is easy to see that
$R(G)$
is an affine algebraic set in $C^{4n}$.See [CS] for details.
If we have a finite volume complete hyperbolic 3-manifold then its holonomy
representation is discrete ,faithful and irreducible.

For each $g \in G$ we define a regular function $\tau_g$ on $R(G)$ by
$\tau_g(\rho)=tr\rho(g)$.
If T is a ring generated by all functions $\tau_g,g\in G$ then
it is finitely generated.See [CS] for a proof.

Also it is shown in [CS] that character space of an irreducible component of
$R(
G)$ is an affine variety.

Thurston([Th] chapter 5) showed that geometric structure of compact manifold M
are determined by holonomy representations of $\pi_1(M)$ near some geometric
 structure
 and its complex dimension is at least $3
\times(\#(generators)-\#(relators)-1)$.

\begin{th}([Th],[CS])Let N be a compact orientable 3-manifold.Let
$\rho:\pi_1(N)
\rightarrow SL_2(C)$ be an irreducible representation such that for each torus
component T of $\partial{N}$,$\rho(im(\pi_1(T) \rightarrow \pi_1(N)))$ not
trivial.
Let R be an irreducible component of $R(\pi_1(N))$ containing $\rho$.Then
 $X_0=character(R)$ has dimension $\geq s-3\chi(N)$,where s is the number of
tori component of $\partial{N}$.
\end{th}
\subsection{Local smooth coordinate chart of Representation space by lengthes
of
finite number of loops}

Let $G$  be a group and $\Re$ be the representation space of $G$ into $SL_2(C)$
,where $\Re$ is the set of equivalence classes of homomorphisms from $G$ into
$SL_2(C)$ and two representations are equivalent if they are conjugate .

Suppose  $\rho$  is a  representation sitting smoothly in
$\Re$.We want to find nice coordinate chart around $\rho$ such that each
coordinate is a translation length of some element in $G$.
\begin{lemma}
Let $\alpha$ be an isometry in $SL_2(C)$. Then
$$
|tr(\alpha)-2|+|tr(\alpha)+2|=2(e^{\frac{l(\alpha)}{2}}+e^{-\frac{l(\alpha)}{
2}}) $$
\end{lemma}
\begin{pf}
 If $\alpha$ is
\[
  \alpha=\left[
           \begin{array}{cc}
            \lambda & 0 \\
              0  &    \lambda^{-1}
           \end{array}  \right]
\]

 Put $x=|\lambda+\lambda^{-1}-2|,y=|\lambda+\lambda^{-1}+2|$.

Then $x^2+y^2=2|\lambda|^2+2\lambda\overline
{\lambda^{-1}}+2\lambda^{-1}\bar{\lambda}
+2|\lambda^{-1}|^2+8$ and
$$x^2+y^2+xy=4(|\lambda|^2+2+|\lambda^{-1}|^2).$$
By using $|\lambda|=e^{\frac{l(\alpha)}{2}}$ claim follows.
\end{pf}

Now let $\acute{G}$ be a subset of $G$ consisting of all elements whose image
under $\rho$ is either hyperbolic or elliptic.
Notice that if $a,b \in{\acute{G}}$ then $a^n,b^n,a^nb^n$ are all in
$\acute{G}$
 for large n.     Also if $a\in G$ and $b\in {\acute{G}}$ then $aba^{-1}$ is in
$\acute{G}$.This implies that marked length spectrum on $G'$ determines
representation up to conjugacy.Choose finite set S of $G$ so that traces of
those
 elements determines representation up to conjugation and trace map from the
neighborhood to ${C^*}^S$ where $C^*=C-\{2,-2\}$
 is an  immersion.This is  possible if  the representation is not
elementary.See
 the arguement after the  example below for justification.

\begin{lemma}(Vogt)
Let G be a free group with three generators.Let $trX_i=x_i,trX_iX_j=y_{ij
},trX_iX_jX_k=z_{ijk}$.Define
$$ P=x_1y_{23}+x_2y_{13}+x_3y_{12}-x_1x_2x_3$$
$$
Q=x_1^2+x_2^2+x_3^2+y_{12}^2+y_{13}^2+y_{23}^2+y_{12}y_{13}y_{23}-x_1x_2y_{12
}\\-x_1x_3y_{13}-x_2x_3y_{23}-4$$
$$\Delta(X_1,X_2,X_3)=P^2-4Q$$
Then $z_{123}$ and $z_{213}$ are roots of the quadratic equation for $z$:
$$z^2-Pz+Q=0$$Given prescribed values $x_i,y_{ij}$ there exists only one
conjugacy class if and only if $\Delta(X_1,X_2,X_3)=0.$
\end{lemma}
See [Ma 1] for the proof.

{\bf{Example}}.Let $G$ be a free group with three generators.I want to show
that
 trace map is not an immersion sometimes.Let  \[ X=\left(
                \begin{array}{cc}
                   \alpha &  0 \\
                  0 & {\alpha}^{-1}
                 \end{array}  \right)  \] and  \[ Y=\left(
                                               \begin{array}{cc}
                                                   \beta & 0 \\
                                                    0    & {\beta}^{-1}
                                                 \end{array} \right)  \]
 Let $Z$ be a third element which does not commute with any of two elements and
$\Delta(X,Y,Z)\neq 0$.
Since $G$ is free neighborhood of this representation is a regular neighborhood
of $SL_2(C) \times SL_2(C) \times SL_2(C) $ and obviously sitting smoothly in
$\Re$.Let tr be a map from this neighborhood to $C^7$ such that
$$tr(\rho)=(trX,trY
,trZ,trXY,trXZ,trYZ,trXYZ)$$Then elementary calculation shows that
\begin{eqnarray*}
dtr_{(X,Y,Z)}(\xi_1,\xi_2,\xi_3)&=&(tr(X\xi_1),tr(Y\xi_2),tr(Z\xi_3),tr(X\xi_2)
\\
   &   & \mbox{} +tr(Y\xi_1),tr(X\xi_3)+tr(Z\xi_1),tr(Y\xi_3)+tr(Z\xi_2),...)
\end{eqnarray*}
where $\xi_i \in sl_2(C).$
  Then  dimension of kernel at $(X,Y,Z)$ is 4.To see this note that there are
two degrees of freedom from $tr(X\xi_1)$ and so on ,so there are six degrees of
freedom from $\xi_1,\xi_2,\xi_3$.

$tr(X\xi_1),
tr(Y\xi_2)$  equal to zero implies $tr(X\xi_2)+tr(Y\xi_1)$ is zero since $X$
and $Y$ commute.But it should satisfy $tr(X\xi_3)+tr(Z\xi_1)=0$ and
$tr(Y\xi_3)+tr(Z\xi_2)=0$ so it takes out two degrees of freedom.
If $\xi_1,\xi_2,\xi_3$ satisfy
$tr(X\xi_1)=tr(Y\xi_2)=tr(Z\xi_3)=tr(X\xi_2)+tr(Y\xi_1)
=tr(X\xi_3)+tr(Z\xi_1)=tr(Y\xi_3)+tr(Z\xi_2)=0$ then the last term vanishes
automatically.Here is the reseason.
If $X_i(t)$ are curves passing throug $X,Y,Z$
respectively and with tangent vectors $\xi_1,\xi_2,\xi_3$ set
$$tr(X_i(t))=x_i(t)
,trX_i(t)X
_j(t)=y_{ij}(t),trX_1(t)X_2(t)X_3(t)=z(t).$$
By lemma 3 $z(t)^2 -P(t)z(t)+Q(t)=0.$If we differentiate this equation at 0 we
get $z'(0)(2z(0)-P(0)=0$ by using $P'(0)=Q'(0)=0$.But $2z(0)-P(0)\neq 0$ since
$\Delta$ is not equal to zero.So tr is not an embedding at this point.

This example has an  obvious geometric meaning.For $tr(X\xi_1)$ to be zero,
diagonal elements of $\xi_1$ should be equal to zero.This implies that we
deform $X$ by
 conjugation.
The same is true for $Y$.Since $X,Y$ have the same axis, conjugating $X$ by
 those elements conjugating $Y$ does not change trace of $X$ so infinitesimal
change of $tr(X)$ is equal to zero. Same argument holds for $Y$.

More explicitly let \[ g_t=\left(
                           \begin{array}{cc}
                            x_1(t) &  x_2(t) \\
                            x_3(t) &  x_4(t)
                           \end{array}
                           \right) \] and  \[ h_t=\left(
                                            \begin{array}{cc}
                                             y_1(t) & y_2(t) \\
                                y_3(t) & y_4(t)
                                             \end{array} \right) \]
such that $g_0=h_0=I$ and $x_2'(0)=-y_2'(0)=x_3'(0)=-y_3'(0)$.
Let \[ Z=\left(
         \begin{array}{cc}
          a &  b  \\
          b &  d  \end{array} \right) \] for simplicity.
Set $X_t=g_tX{g_t}^{-1},Y_t=h_tY{h_t}^{-1},Z_t=g_tZ{g_t}^{-1}$.
Then ${\xi(t)}=(X_t,Y_t,Z_t)$ is a  1-parameter family of defermation which is
not a conjugation defermation.For \[ X'(0)= \left(
                                  \begin{array}{cc}
                                  0 & -x_2'\alpha+x_2'{\alpha}^{-1}\\
                        x_3'\alpha-x_3'{\alpha}^{-1} & 0
                              \end{array} \right) \]
and similarly for $Y'(0)$ there is no conjugation deformation giving tangent
vectors  $X'(0)$ and $Y'(0)$.It is easy to see that $dtr_{(X,Y,Z)}(\xi'(0))=0$
.For if we calculate $\frac{dtr(Y_tZ_t)}{dt}$ at $t=0$,it is equal to
$\beta(x_2'(0)b-x_3'(0)b)+
b(-y_2'(0)\beta+x_2'(0){\beta}^{-1})+b(y_3'(0)\beta-y
_3'(0){\beta}^{-1})+{\beta}^{-1}(-x_2'(0)b+x_3'(0)b)$.

If we look at above example carefully,it shows  if no two elements commute then
trace map
is an immersion,which implies that every tangent vector in the kernel of $dtr$
comes from some conjugation deformation.

Observe the following calculation.
\[  g_t=\left( \begin{array}{cc}
                x_1(t) &  x_2(t) \\
                x_3(t) &   x_4(t)
            \end{array} \right) \]
\[X= \left(
         \begin{array}{cc}
          a &  b  \\
          c &  d  \end{array} \right) \]

Then conjugation of $X$ by $g_t$ give tangent vector at $X$
\[ \left(  \begin{array}{cc}
           x_2'c-x_3'b &  -x_2'a+2x_1'b+x_2'd \\
           x_3'a+2x_4'c-x_3'd   & -x_2'c+x_3'b
     \end{array} \right)  \]
Notice that two entries are linearly independent.So if $g_t,h_t$ give the same
infinitesimal deformation on two different element which do not commute then
they give the same infinitesimal deformation on every element.

By using this fact we can show that if the representation is not elementary
then
 we can show that trace map from a neighborhood of a representation $\rho$
sitting
smoothly          in representation space to $C^N$ for large enough N is a
smooth
immersion.

To see this if $S=\{X_1,\cdots,X_n\}$ is the finite set containing all
generators of the group such that its traces determine representation uniquely
up to conjugacy(which is possible by [CS]) and their images under $\rho$ are
all hyperbolic isometries.By abusing the notation we will identify $X_i$ with
$\rho(X_i)$.
This can be seen that if $S$ is such a set
 possibly
 contaning parabolic or elliptic elements then choose a hyperbolic element
and multiply those elements which are not hyperbolic in $S$ by this element and
its inverse to make them hyperbolic and then throw in this element.Then
using the formula $tr(XY)+tr(XY^{-1})=tr(X)tr(Y)$ we see that all the traces of
old elements can be obtained by new elements.This is a new set consisting
entirely of hyperbolic elements.

Since representation is not elementary there is $X_k$ which does not commute
with $X_1$.Choose any $X_i,i\neq1,k$.If it commutes with $X_1$ replace
$X_1,X_i,
X_k
$ by $X_1X_i,X_1X_k,X_k$.It is easy to check that no two of three does commute.
By    renaming them we get $(X_1,\cdots,X_n)$ such that no two of $X_1,X_2,X_3$
commute
 and $X_i,i\neq 1,2,3$ commutes with at most one of three $X_1,X_2,X_3$.

If $(\xi_1,\cdots,\xi_n)$ is in the kernel of the differential of the trace
map,
there is $g_t$ conjugating  $X_1,X_2,X_3$, giving tagent vectors $\xi_1,\xi_
2,\xi_3$ because no two of three commute.Similary $h_t$ conjugating
$X_1,X_2,X_4
$ giving $\xi_1,\xi_2,\xi_4$ since no two of three commute.But $h_t$ can be
replaced by $g_t$ since they give the same infinitesimal deformation on
$X_1,X_2$ which do not commute.In this way we can see that this tangent
vector
comes from conjugation deformation.

By  lemma 2 the map $f$ from ${C^*}^S$ to $R^{\acute{G}}$ sending each
component $\alpha$ of  an element in ${C^*}^
S$ to $2(e^{\frac{l(\alpha)}{2}}+e^{-\frac{l(\alpha)}{
2}})$ is $C^{\infty}$ since there is no parabolic element in $\acute G$ and
norm
 function is smooth at $C^*$.Furthermore by adding more elements in $S$  we can
make this function an immersion.To see this consider
 the map $g$ from $R^2$ to $R$ such that $g(x,y)=|(x,y)-2|+|(x,y)+2|$ as in the
above lemma.It is easy to see that if $(x,y)\neq  2,-2$ then dimension of the
kernel of $dg$ is one.Suppose $dg(\xi)=0$ at $\alpha$. Then find an element
$\beta$ so that trace($\beta$) is a
polynomial in $\alpha$ and other elements in $S$.Choosing  $\beta$ carefully
it is possible to ensure that $dg(\xi) \neq 0$ at $\beta$.So adding $\beta$ to
$S$ for each such an $\alpha$ we can show the map is an immersion.
\begin{Prop}
Let $\rho$ be an nonelementary representation sitting smoothly inside the
representation space.Then neighborhood of $\rho$,denoted by $N(\rho)$, is
parametrized
 by translation lengthes of finitely many hyperbolic elements in $G$.
\end{Prop}
\begin{pf}
So far we got
$$f\circ tr:N(\rho)\rightarrow {C^*}^S \rightarrow R^{\acute{G}}$$ is an
immersion.But since marked length spectrum on $R^{\acute{G}}$ determines
representation it is one to one also.

Since dimension of
representation space is finite around $\rho$ we can find finite subspace  $R^N$
such that
 projection of the image of the neighborhood to this subspace is injective.
\end{pf}
Specially if the representation is an holonomy of geometrically finite
3-manifold then it is automatically sitting smoothly in representation space
 because Teichmuller
 space  is open in representation space and it satisfies  all the necessary
 condition for above considerations,so small neighborhood of this
representation
 is parametrized smoothly by lengthes of finitely many geodesics.

$$ ACKNOWLEDGEMENT$$
I wish to thank all whom I talked with,among them I want to specially thank
A.Casson as my adviser,C.McMullen,J.P.Otal,F.Bonahon,M.Bourdon,\\
G.Besson for their useful comments.

Inkang Kim \\
University of California,Berkeley\\
Berkeley CA 94720\\
email;inkang@math.berkeley.edu

\end{document}